\documentclass[12pt,a4paper]{article}
\usepackage{amsfonts,amsmath,amssymb,indentfirst}
\newcommand{\er}{\mathbb{R}}

\newtheorem{Teorema}{Theorem}[section]
\newtheorem{Propriedade}[Teorema]{Proposition}
\newtheorem{Rem}[Teorema]{Remark}

\DeclareMathOperator{\sech}{sech}
\newcommand\blfootnote[1]{%
  \begingroup
  \renewcommand\thefootnote{}\footnote{#1}%
  \addtocounter{footnote}{-1}%
  \endgroup
}
\begin{document}
\title{Existence and linearized stability of solitary waves for a quasilinear Benney system}
\author{Jo\~ao-Paulo Dias{\small $^{(a,b)}$}, M\'ario Figueira{\small $^{(a,b)}$} and Filipe Oliveira{\small $^{(a,c)}$}}
\date{{\small(a)CMAF-UL, Av. Prof. Gama Pinto, 1649-003 Lisboa, Portugal}\\
{\small(b) DM-FCUL, Campo Grande, 1749-016 Lisboa, Portugal}\\
{\small(c) FCT-UNL, Monte$\,$ da$\,$ Caparica,$\,$ 2829-516 $\,$Caparica, Portugal}}

\maketitle
\begin{abstract}
\noindent
\blfootnote{\tiny Email addresses: J.P. Dias: dias@ptmat.fc.ul.pt; M.Figueira: figueira@ptmat.fc.ul.pt; F.Oliveira: fso@fct.unl.pt}
We prove the existence of solitary wave solutions to the quasilinear Benney system
\begin{displaymath}
\left\{
 \begin{array}{llll}
 iu_{t}+u_{xx}=a|u|^pu+uv\\
 v_t+f(v)_x=(|u|^2)_x,\\
  \end{array}\right.
\end{displaymath}
where $f(v)=-\gamma v^3$, $-1<p<+\infty$ and $a,\gamma>0$. We establish, in particular, the existence of travelling waves with speed arbitrary large if $p<0$ and
arbitrary close to $0$ if $p>\frac 23$. We also show the existence of standing waves in the case $-1<p\leq \frac 23$, with compact support if $-1<p<0$. 
Finally, we obtain, under certain conditions, the linearized stability of such solutions.
\end{abstract}
\section{Introduction}
\noindent
In the seminal works \cite{B1}, \cite{B2}, D.J. Benney introduced a number of universal models describing the interaction between  short and long waves propagating along a
direction $(Ox)$ in a dispersive media. One of these models is the system
\begin{equation}
\label{Benney1}
 \left\{\begin{array}{llll}
\displaystyle i\frac{\partial u}{\partial t}+\frac{\partial^2 u}{\partial x^2}=m_1|u|^2u+m_2uv\\
\\
\displaystyle\frac{\partial v}{\partial t}+m_3\frac{\partial v}{\partial x}=m_4\frac{\partial}{\partial x}(|u|^2),\quad x\in\er,\,t\geq 0.
        \end{array}\right.
\end{equation}
Here, $m_j$ are real constants, $u=u_{(y)}+iu_{(z)}$ represents, in complex notation, the transverse components $(u_{(y)},u_{(z)})$ of the short wave,
and $v$ the density perturbation induced by the long wave.\\
This model has been successfully applied to several physical contexts, such as the study of the formation 
and annihilation of solitons resulting from the interaction between
Langmuir and ion sound waves in a magnetized plasma, in the case where the perturbation propagates with a speed close to that of sound (\cite{Karpman}, \cite{Yajima}), or
the interaction between Alfv\'en and magneto-acoustic waves in a cold plasma subjected to a strong external magnetic field (\cite{Champeaux}, \cite{Sulem}). In water waves
theory, applications of this model include the interaction between gravity-capillary waves in a two-layer fluid, when the group velocity of the surface waves coincides with the phase velocity of the internal waves 
(see \cite{Funakoshi},\cite{Grimshaw},\cite{Redekopp}. See also \cite{Proment} for an alternative derivation of Benney's equations from the Zakharov formulation of surface gravity waves). Other examples, such as
long-wave short-wave interaction in bubbly liquids (\cite{bubbly}) or optical-microwave interactions in nonlinear mediums (\cite{microwave}) can be given.\\
\\
The mathematical study of system (\ref{Benney1}), namely the well-posedness of the associated Cauchy Problem or the existence and stability of solitary waves, has been extensively conducted over the years by many authors 
(see for instance \cite{angulo}, \cite{Bekiranov}, \cite{Laurencot}, \cite{Ma}, \cite{Tsutsumi1},\cite{Tsutsumi2} and references therein).\\
\\
As pointed out in \cite{B1}, this system is an adequate model in the case where the amplitude of the long wave is considerably smaller than the amplitude of the 
short wave. When both amplitudes are of the same order, the effect of long waves becomes considerably weaker, and,
in this context, (\ref{Benney1}) should be replaced by a system of the form
\begin{equation}
\label{Benney2}
 \left\{\begin{array}{llll}
\displaystyle i\frac{\partial u}{\partial t}+\frac{\partial^2 u}{\partial x^2}=|u|^2u+uv\\
\\
\displaystyle\frac{\partial v}{\partial t}+\frac{\partial }{\partial x}f(v)=\frac{\partial }{\partial x}(|u|^2),
        \end{array}\right.
\end{equation}
where $f$ is a nonlinear polynomial. Contrarely to the linear case $f(v)=mv$, only recently some attention has been given to the mathematical study of these more
general systems. In \cite{Anton}, the case of the Schr\"odinger-Burgers' system ($f(v)=mv^2$) was adressed
in the half-line. The existence and linear stability of shockwave solutions to (\ref{Benney2}) was proved in \cite{Amorim1}. 
By combining methods from dispersive equations and systems of hyperbolic conservation laws, in \cite{DF}, \cite{DFO}, the authors studied the existence of global weak solutions and local strong solutions for the corresponding 
Cauchy problem in the energy space, in the case where $f(v)=av^2-bv^3$, $a,b>0$ 
(see also \cite{nAD}, \cite{Amorim2}, \cite{nDFF1}, \cite{nDFF2}, \cite{nDF} and \cite{nDFO} for related results concerning similar systems).
\\
\\
Also recently, B\'egout and D\'\i az (\cite{BD1}, \cite{BD2}) considered nonlinear Schr\"odinger equations with an ``absorbing'' singular potential of the form $|u|^p$, $p<0$ such as the homogenous equation
$$iu_t+\Delta u=\alpha|u|^{p}u,\qquad -1<p<0.$$
Nonlinear Schr\"odinger equations with singular potentials arise in a large variety of contexts (see e.g. \cite{LeMesurier},\cite{Stollmann}).
The authors proved in particular that under some circumstances such equations admit standing wave solutions of the form $u(x,t)=\phi(x)e^{i\beta t}$ with compact support,
under the fundamental condition $-1<p<0$. Such localization of solutions is well-known not to exist for ordinary Schr\"odinger equations and seem to be
a special feature of singular potentials of this type.\\
\\
With these motivations, in the present work, we are concerned with the existence and behaviour of solitary waves for quasilinear Benney systems of the type
\begin{equation}
\label{equacaoinicial}
 \left\{\begin{array}{llll}
\displaystyle i\frac{\partial u}{\partial t}+\frac{\partial^2 u}{\partial x^2}=m_1|u|^pu+uv\\
\\
\displaystyle\frac{\partial v}{\partial t}+\frac{\partial }{\partial x}f(v)=\frac{\partial }{\partial x}(|u|^2),
        \end{array}\right.
\end{equation}
where $f(v)=m_2v^3$ and $-1<p<+\infty$.\\
\\
\\
The rest of this paper is organized as follows:\\
\\
In Sections 2 and 3 we establish the existence of a two-parameter family of solitary-wave solutions to (\ref{equacaoinicial}) of the form
\begin{equation}
 \label{forma}
(u(x,t),v(x,t))=(e^{iwt}e^{i\frac c2(x-ct)}\phi(x-ct),{\psi}(x-ct)),
\end{equation}
where $\phi$ and $-\psi$ are non-negative radially decreasing functions vanishing at infinity. This result relies on the derivation of sharp estimates for the Lagrange multiplier associated to
a variational minimization problem. These estimates allow us also to exhibit solitary waves with positive speed $c$ arbitrary large in the case $-1<p<0$ and arbitrary close to $0$ for 
$p>\frac 23$.\\
\\
When $0\leq p\leq \frac 23$, we prove, in Section 4, the existence of standing-wave solutions ($c=0$) of the form
$$(u(x,t),v(x,t))=(e^{iwt}\phi(x),{\psi}(x))$$
by applying a result due to Berestycki and Lions (\cite{BerLions}). We also establish the existence of standing waves with compact support. The condition for the existence of such localized solutions is $-1<p<0$, related in particular to the convergence of a singular integral of the type
$\displaystyle \int_0^{a}\frac{dx}{x^{1+\frac p2}}$. Although we use totally different methods, this is, as mentionned above, the exact same condition used in \cite{BD1}, \cite{BD2} to derive solutions with compact support.\\
\\
Finally, in Section 5, after establishing the global well-posedness of a non-autonomous system consisting of the linearization of (\ref{equacaoinicial}) around a solitary wave, we prove, in the spirit of \cite{GR}, the linearized stability of solitary wave solutions in the case $p>-\frac 23$, with $c=0$ if $p<0$ (and without restrictions on the speed $c$
if $p>0$).\\
Our results are synthesized in the following table: 
\begin{center}
 {\bf Existence and stability of solitary-wave solutions to (\ref{equacaoinicial})}
\end{center}
\begin{table}[h]
\begin{tabular}{|c|c|l|l|l|}
\hline
 &Regularity  &Speed  &Further Properties    \\ \hline
 $-1<p<0$&$H^1\times H^1$  &$\left.\begin{array}{llll}
                              c>0,\\
                             \textrm{arbitrarely large}
                            \end{array}\right.$ &    \\ \hline
$-1<p\leq -\frac 23$&$\left.\begin{array}{llll}
                              (C^2\cap W^{2,\infty })\\
                              \times(C^1\cap W^{1,\infty })
                           
                            \end{array}\right.$   &$\left.\begin{array}{llll}
                              c=0\\
                             \textrm{}
                            \end{array}\right.$&$\left.\begin{array}{llll}
                             \textrm{Compactly}\\
                             \textrm{supported }
                            \end{array}\right.$   \\ \hline
 $-\frac 23<p\leq \frac 23$&$\left.\begin{array}{llll}
                              (C^3\cap W^{3,\infty })\\
                              \times(C^2\cap W^{2,\infty })
                           
                            \end{array}\right.$ 
 &$\left.\begin{array}{llll}
                              c=0\\
                             \textrm{}
                            \end{array}\right.$ &$\left.\begin{array}{llll}
                              \textrm{Linearly stable;}\\
                             \textrm{Compactly}\\
                             \textrm{supported if $p<0$ }
                            \end{array}\right.$  \\ \hline
 $p>\frac 23$&$H^{\infty}\times H^{\infty}$  &$\left.\begin{array}{llll}
                              c>0,\\
                             \textrm{arbitrarely small}
                            \end{array}\right.$   & $\textrm{ Linearly stable}$  \\ \hline
\end{tabular}
\end{table}
\section{Existence of Solitary waves for $-1<p<0$}
\noindent
We consider the system
\begin{equation}\label{sistemainicial}
\left\{
 \begin{array}{llll}
 iu_{t}+u_{xx}=a|u|^pu+uv\\
 v_t+f(v)_x=(|u|^2)_x,\\
  \end{array}\right.
\end{equation}
where $f(v)=-\gamma v^3$, $-1<p<0$, $\gamma>0$ and $a>0$.\\
\\
We look for solutions of the form 
$$(u(x,t),v(x,t))=(e^{iwt}e^{i\frac c2(x-ct)}\phi(x-ct),{\psi}(x-ct)),$$ 
with $\phi$ and $\psi$ real-valued and vanishing at infinity. 
We obtain the system
\begin{equation}\label{sistemaest}
\left\{
 \begin{array}{rrrr}
-\phi''+c^*\phi&=&-\phi\psi-a|\phi|^p\phi\\
 c\psi&=&-\phi^2+f(\psi),\\
  \end{array}\right.
\end{equation}
where $c^*=w-\frac{c^2}4$.\\
\\
By showing the existence of solutions to (\ref{sistemaest}), we will prove the following theorem, describing a two-parameter family of soutions to (\ref{sistemainicial}):
\begin{Teorema}
\label{principal}
Let $\displaystyle\frac 13<\alpha< 1$.\\
\\
There exists $\mu_0=\mu(\alpha)>0$ such that for all $\mu>\mu_0$, the system (\ref{sistemainicial}) has non-trivial solutions of the form
\begin{displaymath}
 \left\{
 \begin{array}{llll}
  u(x,t)=e^{iwt}e^{i\frac c2(x-ct)}\phi_{\mu,\alpha}(x-ct),\\
  \\
  v(x,t)={\psi}_{\mu,\alpha}(x-ct)
 \end{array}\right.
\end{displaymath}
where $\phi_{\mu,\alpha}$ and $-\psi_{\mu,\alpha}$ are non-negative radially decreasing $H^1$ functions such that
$$\|\phi_{\mu,\alpha}\|_{H^1}^2+\|\psi_{\mu,\alpha}\|_2^2\geq \mu^{\frac 14(1-\alpha)}.$$
Furthermore,
$$c=c(\mu,\alpha)\approx_{\mu\to+\infty}\mu^{\frac 12(3-\alpha)}.$$
\end{Teorema}

\bigskip

\noindent
{\bf The minimization problem}\\
\\
For $u\in  H^1(\er)\cap L^{p+2}(\er)$, $-1<p<0$, and $v\in L^2(\er)\cap L^4(\er)$, let
$$\tau(u,v)=\frac {2a}{p+2}\int |u|^{p+2}+\int vu^2+\frac {\gamma}4\int v^4.$$
Also, for $d,\mu>0$,  let
\begin{multline*}
 X_{\mu,d}=\{(u,v)\in H^1(\er)\cap L^{p+2}(\er)\times (L^2(\er)\times L^4(\er)) \,:\\
\,{N_d}(u,v)=\|u\|_2^2+\|u'\|_2^2+d\|v\|_2^2=\mu\}
\end{multline*}
and
\begin{equation}
 \label{mp}
\mathcal{I}(\mu,d)=inf \{\tau(u,v)\,:\,(u,v)\in X_{\mu,d}\}.
\end{equation}
\noindent
If $(u,v)$ is a minimizer, then there exists a Lagrange multiplier $\lambda$ such that
$\nabla \tau=\lambda \nabla N_d$, that is

\begin{equation}\label{eulerlagrange}
\left\{
 \begin{array}{lllll}
2a|u|^pu+2vu&=&\lambda(-2u''+2u)\\
u^2+\gamma v^3&=&2\lambda dv
  \end{array}\right.
\end{equation}
and
\begin{equation}\label{eulerlagrange2}
\left\{
 \begin{array}{lllllllll}
\lambda u''-\lambda u&=&-uv-a|u|^pu\\
-2d\lambda v&=&-u^2+f(v).
  \end{array}\right.
\end{equation}
If $\lambda<0$, the change of variable $x'=x\sqrt{-\lambda}$ leads to a solution
\begin{equation}
 \label{scale}
(\phi(x),\psi(x))=\left(u\left(\sqrt{-\lambda}x\right),v\left(\sqrt{-\lambda}x\right)\right)
\end{equation}
of system (\ref{sistemaest}) for
\begin{equation}
\label{velocidades}
c*=-\lambda\quad\textrm{ and }\quad c=-2\lambda d.
\end{equation}

\bigskip 

\begin{Propriedade}
\label{neinfinito}
 For $\mu,d>0$, $\mathcal{I}(\mu,d)>-\infty$.
\end{Propriedade}
{\bf Proof:}\\
We only have to notice that for $(u,v)\in X_{\mu,d}$, $$\tau(u,v)\geq - \int |v|u^2\geq -\|v\|_2\|u\|_4^2\geq -C\|v\|_2\|u'\|_2^{\frac 12}\|u\|_2^{\frac 32}
\geq -C\frac{\mu^{\frac 32}}{d^{\frac 12}},$$
by the Gagliardo-Nirenberg inequality ($C>0$).\hfill$\blacksquare$

\bigskip

\begin{Propriedade}
\label{estneg}
For $\mu,d>0$, $$\displaystyle \mathcal{I}(\mu,d)\leq - \frac{3}{8\sqrt{\pi}}\frac{\mu^{\frac 32}}{d^{\frac 12}}+C\left(\mu^{1+\frac p2}+
{\gamma}\frac{\mu^2}{d^2}\right),$$ where $C$ is a positive constant.\\
In particular, for $\displaystyle\frac 13<\alpha<1$, $d=\mu^{\alpha}$ and $\mu$ large enough, $\mathcal{I}(\mu,d)<0$.
\end{Propriedade}
{\bf Proof: }\\
For $B>0$, we consider the following functions
$$u(x)=\frac{B}{1+x^2}$$
and
$$v(x)=-\frac 1{\sqrt{d}}u(x).$$
A simple computation shows that
$$\|u\|_2^2+\|u'\|_2^2+d\|v\|_2^2=B^2\pi,$$
hence, by taking $\displaystyle B=\sqrt{\frac\mu\pi}$, $(u,v)\in X_{\mu,d}$.\\
\\
Furthermore, $$\displaystyle \int vu^2=-\frac{B^3}{\sqrt{d}}\int\left(\frac 1{1+x^2}\right)^3=
-\frac{\frac {3\pi}8}{\pi^{\frac 32}}\frac{\mu^{\frac 32}}{d^{\frac 12}}=-\frac{3}{8\sqrt{\pi}}\frac{\mu^{\frac 32}}{d^{\frac 12}},$$ hence
$$\tau(u,v)=\frac {2a}{p+2}\int |u|^{p+2}+\int vu^2+\frac {\gamma}4\int v^4
\leq -\frac{3}{8\sqrt{\pi}}\frac{\mu^{\frac 32}}{d^{\frac 12}}+C\left(\mu^{1+\frac p2}+{\gamma}\frac{\mu^2}{d^2}\right),$$
where $C\displaystyle >0$.\hfill$\blacksquare$

\noindent

\begin{Propriedade}\label{lemaradial}
Let $\mu,d>0$ and $(u,v)\in X_{\mu,d}$.\\
There exists $\tilde{u}$ non-negative and $\tilde{v}$ non-positive, $\tilde{u}$ and $|\tilde{v}|$ radially decreasing,
 such that $\tau(\tilde{u},\tilde{v})\leq \tau(u,v)$ and $(\tilde{u},\tilde{v})\in X_{\mu,d}$.  
\end{Propriedade}
{\bf Proof:}\\
Let $u_*=|u|^*$ and $v_*=-|v|^*$, where $f^*$ denotes the Schwarz symmetrization of $f$.\\
On one hand, $$\tau(|u|,-|v|)=\frac {2a}{p+2}\int |u|^{p+2}-\int |v|u^2+\frac {\gamma}4\int v^4\leq \tau(u,v).$$
Furthermore, since for $r\geq 1$, $\displaystyle \int (f^*)^r=\int f^r$ for every positive function $f$ in $L^r(\er)$ and $\displaystyle \int |u|^2|v|\leq \int (|u|^*)^2|v|^*$,
$$\tau(u_*,v_*)\leq \tau(u,v).$$
By the  Polya-Szego inequality, $\displaystyle \int ((u_*)')^2\leq \int (u')^2$, hence $$N_d(u_*,v_*)\leq N_d(u,v)=\mu.$$
\\
If $N_d(u_*,v_*)=\mu$, we put $(\tilde{u},\tilde{v})=(u_*,v_*)$.\\
If $N_d(u_*,v_*)<\mu$ we set, for $k>0$,
\begin{equation}
\displaystyle \tilde{u}(x)=k^{\frac 1{4p}}u_*\left(\frac x{k^{\frac{p+2}{-4p}}}\right) \textrm{ and } \displaystyle \tilde{v}(x)=k^{\frac 14}v_*(kx).
\end{equation}
\\
Since $\displaystyle \int |\tilde u|^2=k^{-\frac 14}\int u_*^2$ and $\displaystyle \int |\tilde v|^2=k^{-\frac 12}\int v_*^2$ and at least one of these quantities is different
from $0$, there exists $0<k<1$ such that $$N_d(\tilde{u},\tilde{v})=\mu.$$
Furthermore,
$$\int \tilde{v}^4=\int v_*^4,$$
$$\int \tilde{u}^{p+2}=\int u_*^{p+2}$$
and
$$\displaystyle \int \tilde{u}^2\tilde{v}=k^{\frac 1{2p}+\frac 14}\int u_*^2\left(\frac x{k^{\frac{p+2}{-4p}}}\right)v_*(kx)< k^{\frac 1{2p}+\frac 14}\int u_*^2\left(\frac x{k^{\frac{p+2}{-4p}}}\right)v_*\left(\frac x{k^{\frac{p+2}{-4p}}}\right)$$
since 
\begin{equation}
 \label{kapa}
\displaystyle |kx|< \left|\frac x{k^{\frac{p+2}{-4p}}}\right|
\end{equation}
for $x\neq 0$ and  $-v_*$ is non-negative and radially decreasing. Finally,
$$\displaystyle \int \tilde{u}^2\tilde{v}\leq \left(k^{\frac 1{2p}+\frac 14-{\frac{p+2}{4p}}}\right)\int u_*^2v_*=\int u_*^2v_*$$
and $\tau(\tilde{u},\tilde{v})< \tau(u_*,v_*)\leq \tau(u,v)$, which completes the proof.\hfill$\blacksquare$\\
\\
\begin{Propriedade}
\label{minimizador}Let $\mu,d>0$. There exists a solution $(u,v)$ for the minimization problem (\ref{mp}), with $u$ and $-v$ non-negative and radially decreasing.
\end{Propriedade}
{\bf Proof:}\\
Let $(u_n,v_n)$ a minimizing sequence in $(H^1_{rd}(\er)\cap L^{p+2}(\er))\times L^2_{rd}(\er)\cap L^4(\er)$.\\
By the compacity of the injection $H^1_{rd}(\er)\hookrightarrow L^r(\er)$, $r>2$, there exists a subsequence still denoted $u_n$ such that
\begin{itemize}
\item $u_n\to u$ in $L^4(\er)$;
\item $u_n\rightharpoonup u$ in $H^1(\er)$ weak;
\item $u_n\to u$ almost everywhere (in particular, $u$ is radial decreasing).\\
\end{itemize}
Also, since $\displaystyle \|v_n\|_2^2\leq \frac{\mu}{d}$ is bounded, we can extract a subsequence still denoted $v_n$ such that 
 $v_n\rightharpoonup v$ in $L^2(\er)$ weak.\\
 \\
Hence, since $u_n^2\to u^2$ in $L^2$ strong and $v_n\rightharpoonup v$ in $L^2$ weak, $$\int v_nu_n^2\to \int u^2v.$$
The sequence $$\frac{\gamma}4\int v_n^4=\tau(u_n,v_n)-\int v_nu_n^2-\frac{2a}{p+2}\int |u|^{p+2}$$
is thus bounded, and we can extract a subsequence still denoted $v_n$ such that 
$v_n\rightharpoonup v$ in $L^4$ weak.\\
\\
Since $\displaystyle \int v^4\leq \liminf \int v_n^4$ and $\displaystyle\int |u|^{p+2}\leq \liminf \int |u_n|^{p+2}$,
$$\tau(u,v)\leq \liminf \tau(u_n,v_n)=\mathcal{I}(u,v).$$
Now, if $\|u\|_2^2+\|u'\|_2^2+d\|v\|_2^2<\mu$, the construction made in the proof of Proposition \ref{lemaradial} shows that there exists
$(\tilde{u},\tilde{v})\in X_{\mu,d}$ such that $$\tau(\tilde{u},\tilde{v})<\tau(u,v).$$
Finally, $$\mathcal{I}(\mu,d)\leq \tau(\tilde{u},\tilde{v})<\tau(u,v)\leq \liminf \tau(u_n,v_n)=\mathcal{I}(\mu,d),$$
which is absurd, hence $(\tilde{u},\tilde{v})\in X_{\mu,d}$ is a minimizer.\\
\\
Note that, since $\mathcal{I}(\mu,d)=\tau(u,v)$, we have in fact that $\int |u|^{p+2}=\liminf \int |u_n|^{p+2}$, $\int v^4=\liminf \int v_n^4$, 
 $\int v^2=\liminf \int v_n^2$ and $\int u^2+u'^2=\liminf \int u_n^2+u_n'^2$, hence $u_n\to u$ in $L^{p+2}(\er)\cap H^1(\er)$ strong and $v_n\to v$ in 
 $L^2(\er)\cap L^4(\er)$ strong.\\
 By choosing a new subsequence, $v_n\to v$ almost everywhere, hence $-v$ is non-negative and radially decreasing.\hfill$\blacksquare$
\\
\\
If $(u,v)\in X_{\mu,d}$ is a solution to the minimization problem, $u,-v\geq 0$, there exists a Lagrange multiplier $\lambda\in \er$ such that 
\begin{equation}\label{eulerlagrange3}
\left\{
 \begin{array}{rrrr}
\lambda u''-\lambda u&=&-uv-a|u|^pu\\
-2d\lambda v&=&-u^2+f(v).
  \end{array}\right.
\end{equation}
The next result states the assymptotic behaviour of $\lambda$:
\begin{Propriedade}
\label{asym}
Let $\displaystyle\frac 13<\alpha < 1$ and $d=\mu^{\alpha}$. There exists positive constants $M_1,M_2$ such that for $\mu$ large enough, 
$$M_1 \left(\frac{\mu}{d}\right)^{\frac 32}\leq -\lambda \leq M_2\left(\frac{\mu}{d}\right)^{\frac 32}.$$
\end{Propriedade}
{\bf Proof:}\\
Multiplying the equations in (\ref{eulerlagrange3}) respectively by $\phi$ and $\psi$,
\begin{equation}
\label{estimativabase}
2\lambda d=3\int u^2v+2a\int u^{p+2}+\gamma \int v^4.
\end{equation}
In particular,
$$-2\lambda d=-3\int u^2v-2a\int u^{p+2}-\gamma\int v^4$$
\begin{equation}
\label{estsup}
\leq 3\left(\int u^4\right)^{\frac 12}\left(\int v^2\right)^{\frac 12}\leq 3C\mu\sqrt{\frac {\mu}d},
\end{equation}
$C>0$, by the Gagliardo Nirenberg inequality.\\
This proves the second inequality by choosing $\displaystyle M_2=\frac {3C}2$.\\
\\
Now, since $$2a\int u^{p+2}=(p+2)\tau(u,v)-(p+2)\int u^2v-(p+2)\frac{\gamma}4\int v^4,$$
we obtain by (\ref{estimativabase}) that
\begin{equation}
 \label{estinf1}
2\lambda d=(1-p)\int u^2v+(p+2)\tau(u,v)+\frac{\gamma}4\left(2-p\right) \int v^4.
\end{equation}
Since 
\begin{equation}
 \label{estu4}
\frac{\gamma}4\int v^4=\tau(u,v)-\frac {2a}{p+2}\int u^{p+2}-\int u^2v\leq \tau(u,v)-\int u^2v,
\end{equation}
$$2\lambda d\leq 4\tau(u,v)-\int u^2v\leq 4\tau(u,v)+\left(\int v^2\right)^{\frac 12}\left(\int u^4\right)^{\frac 12}$$
$$\leq 4\tau(u,v)
+C_0^{\frac 12}\frac {\mu^{\frac 32}}{d^{\frac12}},$$
where $C_0$ is the smaller constant for the Gagliardo-Nirenberg inequality $\|u\|^4_4\leq C_0\|u'\|_2\|u\|_2^3.$\\
\\
By Proposition $\ref{estneg}$,
\begin{equation}
\label{estinffinal}
2\lambda d\leq \left(C_0^{\frac 12}-\frac{3}{2\sqrt{\pi}}\right)\frac {\mu^{\frac 32}}{d^{\frac12}}+4C\left(\mu^{1+\frac p2}+
{\gamma}\frac{\mu^2}{d^2}\right),
\end{equation}
where $C>0$. Also, one can choose $C_0=\frac1{\sqrt{3}}$. Indeed, it is known that the sharp constant in the Gagliardo-Nirenberg inequality is given by $C_0=\frac 4{\sqrt{3}\|Q\|_2^2}$, where $Q(x)=\sqrt{2}\sech(x)$ is
the positive radial solution of $Q''+Q^3=Q$: $\|Q\|_2^2=4$ (see for instance \cite{GN1}, \cite{GN2}).
\\
\\
Now, taking $d=\mu^{\alpha}$ and putting $\displaystyle \epsilon=\frac{3}{2\sqrt{\pi}}-\frac 1{3^{\frac 14}}>0$,
$$c:=-2\lambda d\geq \epsilon \frac {\mu^{\frac 32}}{d^{\frac12}}-C'\left(\mu^{1+\frac p2}+
\frac{\gamma}{2}\frac{\mu^2}{d^2}\right)=\epsilon \mu^{\frac 32-\frac \alpha 2}-C'\mu^{1+\frac p2}-C'\frac{\gamma}2\mu^{2(1-\alpha)}$$
$$\geq \frac{\epsilon}{2}\mu^{\frac 32-\frac \alpha 2}$$
for $\mu$ large enough, since for $1\geq \alpha>\frac 13$, we have $1+\frac{p}2<\frac 12(3-\alpha)$ and $2(1-\alpha)<\frac 12(3-\alpha)$.\\
The proof is now complete by taking $\displaystyle M_1=\frac{\epsilon}4$.\hfill $\blacksquare$
\\
\\
{\bf End of the proof of Theorem \ref{principal}:}\\
\\
In particular, from Proposition \ref{asym}, $\lambda<0$. By the change of variables (\ref{scale}), we obtain from a minimizer $(u,v)\in X_{\mu,d}$ a solution $(\phi_\mu,\psi_\mu)$ of system (\ref{sistemaest}).\\
Note that
$$\mu=\|u\|_2^2+\|u'\|_2^2+d\|v\|_2^2=\left\|\phi_{\mu,\alpha}\left(\frac{\cdot}{\sqrt{-\lambda}}\right)\right\|_2^2+\left\|\phi_{\mu,\alpha}'\left(\frac \cdot{\sqrt{-\lambda}}\right)
\right\|_2^2+d\left\|\psi_{\mu,\alpha}\left(\frac \cdot{\sqrt{-\lambda}}\right)\right\|_2^2$$
and
$$\mu=\sqrt{-\lambda}\|\phi_{\mu,\alpha}\|_2^2+\frac 1{\sqrt{-\lambda}}\|\phi_{\mu,\alpha}'\|_2^2+d\sqrt{-\lambda}\|\psi_{\mu,\alpha}\|_2^2.$$
Hence,
\begin{equation}
\label{equivalente}
\mu\approx \mu^{\frac 34(1-\alpha)}\|\phi_{\mu,\alpha}\|_2^2+ \mu^{\frac 34(\alpha-1)}\|\phi_{\mu,\alpha}'\|_2^2+\mu^{\frac 14(3+\alpha)}\|\psi_{\mu,\alpha}\|_2^2
\end{equation}
$$\leq\mu^{\frac 14(3+\alpha)}\left(\|\phi_\mu\|_{H^1}^2+\|\psi_\mu\|_2^2\right)$$
and
$$\|\phi_{\mu,\alpha}\|_{H^1}^2+\|\psi_{\mu,\alpha}\|_2^2\geq C\mu^{\frac 14(1-\alpha)},\quad C>0.$$
$\,$\hfill$\blacksquare$

\section{Existence of Solitary waves for $p>\frac 23$}
\noindent
In the case of $p>\frac 23$, we prove the following result:
\begin{Teorema}
\label{principal2}
Let $\displaystyle 1-p<\alpha<\frac 13$.\\
\\
There exists $\mu_0=\mu(\alpha)>0$ such that for all $0<\mu<\mu_0$, the system (\ref{sistemainicial}) has non-trivial solutions of the form
\begin{displaymath}
 \left\{
 \begin{array}{llll}
  u(x,t)=e^{iwt}e^{i\frac c2(x-ct)}\phi_{\mu,\alpha}(x-ct),\\
  \\
  v(x,t)={\psi}_{\mu,\alpha}(x-ct)
 \end{array}\right.
\end{displaymath}
where $\phi_{\mu,\alpha}$ and $-\psi_{\mu,\alpha}$ are non-negative radially decreasing smooth functions such that,  for $0\displaystyle <\alpha<\frac 13$,
\begin{equation}
\label{final}
\|\phi_{\mu,\alpha}\|_{H^1}^2+\|\psi_{\mu,\alpha}\|_2^2\leq C\mu^{\frac 14(1-3\alpha)},\quad C>0.
\end{equation}
Furthermore,
\begin{equation}
\label{velocidade}
c=c(\mu,\alpha)\approx_{\mu\to 0^+}\mu^{\frac 12(3-\alpha)}.
\end{equation}
\end{Teorema}
\noindent
{\bf Proof:}\\
We begin by noticing that Propositions \ref{neinfinito} and \ref{estneg} hold for $\displaystyle p>\frac 23$. Furthermore, estimate (\ref{kapa}) holds for $p<0$ and for 
$\displaystyle p>\frac 23$.\\
Hence, the conclusions in Propositions \ref{lemaradial} and \ref{minimizador} can be drawn also in this case.\\
Finally, estimate (\ref{estsup})
$$-\lambda d\leq \frac{3C}2\frac{\mu^{\frac 32}}{d^{\frac 12}}$$
remains valid for all $p$, and, for $2-p\geq 0$, estimate (\ref{estinffinal})
$$2\lambda d\leq \left(C_0^{\frac 12}-\frac{3}{2\sqrt{\pi}}\right)\frac {\mu^{\frac 32}}{d^{\frac12}}+4C\left(\mu^{1+\frac p2}+
{\gamma}\frac{\mu^2}{d^2}\right),$$
can be derived in the exact same way as in the case $p<0$. 
On the other hand, if $p> 2$, we get from (\ref{estinf1}) and (\ref{estu4}) that 
$$2\lambda d=4\tau(u,v)-\int u^2v+2a\frac{p-2}{p+2}\int u^{p+2}.$$
By Proposition \ref{estneg},
$$2\lambda d\leq \left(C_0^{\frac 12}-\frac{3}{2\sqrt{\pi}}\right)\frac {\mu^{\frac 32}}{d^{\frac12}}+4C\left(\mu^{1+\frac p2}+
{\gamma}\frac{\mu^2}{d^2}\right)+2a\frac{p-2}{p+2}\int u^{p+2}.$$
Using the Gagliardo-Nirenberg inequality $\|u\|_{p+2}\leq C\|u\|^{\frac p{2p+4}}\|u'\|^{\frac{p+4}{2p+4}}$, we obtain $\displaystyle \int u^{p+2}\leq C\mu^{1+\frac p2},$
hence, in all cases, 
$$c=-2\lambda d\geq \epsilon \mu^{\frac 32-\frac \alpha 2}-C_1\mu^{1+\frac p2}-C_2\frac{\gamma}2\mu^{2(1-\alpha)},$$
where $\epsilon, C_1,$ and $C_2$ are positive constants and $d=\mu^{\alpha}$.\\
\\
Taking $\displaystyle 1-p<\alpha<\frac 13$, $\displaystyle 1+\frac p2>\frac 32-\frac{\alpha}2$ and $\displaystyle 2(1-\alpha)>\frac 32-\frac{\alpha}2$, hence there exists
$\mu_0>0$ such that for all $0<\mu<\mu_0$, 
$$c\geq \frac {\epsilon}2\mu^{\frac 32-\frac{\alpha}2},$$
which, with estimate $(\ref{estsup})$, completes the proof of (\ref{velocidade}).\\
Finally, estimate (\ref{final}) follows from (\ref{equivalente}).
\begin{Rem}
 \label{regularidade}
 In whats concerns the regularity of $\phi$ and $\psi$, note that the monotony of $\phi$ and $\psi$ garantee, via Lebesgue's Theorem, that $\phi'$ and $\psi'$ exist
 almost everywhere. Differentiating the second equation in (\ref{sistemaest}) then yields
 $$\psi'(c+3\gamma\psi^2)=-2\phi\phi'.$$
 Since $\phi\in H^1(\er)$,
$$\int (\psi')^2\leq \frac 4{c^2} \int \phi^2\phi'^2\leq \frac 4{c^2} \|\phi\|_{\infty}^2\|\phi'\|_2^2<+\infty$$
and $\psi \in H^1(\er)$.\\
Now, in the case where $p\geq 0$, the first equation in (\ref{sistemaest}) shows that $\phi''\in L^2(\er)$, that is, $\phi\in H^2(\er)$. And again, by differentiating the second equation,
$$\psi''(c+3\gamma \psi^2)=-2(\phi')^2-2\phi\phi''-6\gamma \psi(\psi')^2,$$
and we easily get that in fact $\psi\in H^2(\er)$. A bootstrap argument then shows that in this case $\phi,\psi\in H^{\infty}(\er)$.
\end{Rem}
\begin{Rem}
\label{remu}
For $p\geq 0$ and $c\geq0$, let $ (\phi,\psi)$ be $C^2(\er)\cap W^{2,\infty}$ solutions of (\ref{sistemaest}) with $\phi\geq 0$ and $\psi\leq 0$. Then $\phi^p\in C^2(\er)\cap W^{2,\infty}$.\\
Indeed, in a neighbourhood of a point $x$ such that
$\phi(x)=0$, from the second equation in (\ref{sistemaest}), $\psi\sim\phi^2$ if $c>0$ ($\psi\sim\phi^{\frac 23}$ if $c=0$). Hence, we derive from the first
equation in (\ref{sistemaest}) that $\phi''\sim \phi$. Noticing that $\phi$ is non-negative and non-increasing, $\phi(x)=0$ implies that $\phi(y)=0$ if $y>x$ and, in
particular, $\phi'(x)=0$. Writing $\phi'(y)=\int_x^y\phi''(t)dt$ then shows that $\phi'\sim \phi$. Finally, $\phi^{p-1}\phi'$, $\phi^{p-2}(\phi')^2$ and
$\phi^{p-1}\phi''$ vanish at $x$, which gives the desired result.
\end{Rem}

\section{Existence of standing waves for $-1< p\leq \frac 23$}
\noindent
In this section we show the existence of smooth non-trivial standing wave solutions to (\ref{equacaoinicial}).
More precisely:
\begin{Propriedade}
 Let $\displaystyle -1<p\leq \frac 23$ and $\displaystyle \gamma,a,\omega>0$, with $\displaystyle \gamma^{-\frac 13}>a$ if $p=\frac 23$.\\
 Then (\ref{equacaoinicial}) admits non trivial solutions of the form 
 $$(u(x,t),v(x,t))=(e^{iwt}\phi(x),\psi(x)),$$
 where $\phi\in C^2(\er)\cap W^{2,\infty}(\er)$ and $-\psi=\left(\frac{\phi^2}{\gamma}\right)^{\frac 13}\in C^1(\er)\cap W^{1,\infty}(\er)$ are non-negative,
 radially descreasing functions.\\
 Moreover:
 \begin{itemize}
  \item if $p>-\frac 23$, $\phi\in C^3(\er)\cap W^{3,\infty}(\er)$ and $\psi\in C^2(\er)\cap W^{2,\infty}(\er)$;
 \item if $-1<p<0$, $\phi$ and $\psi$ are compactly supported.
\end{itemize}
\end{Propriedade}
\noindent{\bf Proof:}\\
Let us consider the system (\ref{sistemaest}) with $c=0$ and $\phi\geq 0$:
\begin{equation}\label{sistemaest2}
\left\{
 \begin{array}{cccc}
-\phi''+w\phi&=&-\phi\psi-a\phi^{p+1}\\
 \phi^2&=&-\gamma \psi^3.\\
  \end{array}\right.
\end{equation}
From the second equation, we obtain $\psi=-\left(\frac{\phi^2}{\gamma}\right)^{\frac 13}$.\\
Replacing in the first equation leads to
\begin{equation}
 \label{equacao1}
 \phi''=a\phi^{p+1}+\omega \phi-\gamma^{-\frac 13}\phi^{\frac 53}.
\end{equation}
We first analyse the case $-1<p< 0$.\\
By multiplying (\ref{equacao1}) by $\phi'$ and integrating, we deduce, for a solution verifying $\phi'(\xi)=0$ in all points $\xi$ such that $\phi(\xi)=0$, that
\begin{equation}
\label{equacao2}
\phi'^2=\frac{2a}{p+2}\phi^{p+2}+w\phi^2-\frac 34\gamma^{-\frac 13}\phi^{\frac 83}:=h(\phi).
\end{equation}
Now, taking $\phi_0>0$ such that $h(\phi_0)=0$ and $h(\phi)\neq 0$ for $\phi\in]0,\phi_0[$, we can derive from (\ref{equacao2}) the existence of  a solution $\phi\in C^2(\er)$ to (\ref{equacao1}) with compact support,
non-negative, radially decreasing, such that $\max \phi=\phi(0)=\phi_0$ and $supp(\phi)=[-x_0,x_0]$, $x_0=\int_0^{\phi_0}(h(\phi))^{-\frac 12}d\phi$
(note that this integral
 is finite for $p<0$).\\
Moreover, if $\displaystyle -\frac 23<p<0$, we can easily establish, from (\ref{equacao1}) and (\ref{equacao2}), that $\phi\in C^3(\er)$ and $\psi\in C^2(\er)$, with the same support.\\
\\
We now turn to the case $\displaystyle 0\leq p\leq \frac 23$, with $\displaystyle \gamma^{-\frac 13}>a$ if $\displaystyle p=\frac 23$. Equation (\ref{equacao1}) can be written as
\begin{equation}
 \label{nova21}
 -\phi''=g(\phi):=-a\phi^{p+1}-w\phi+\gamma^{-\frac 13}\phi^{\frac 53}.
\end{equation}
We have $g\in C^1(\er)$, $g(0)=0$ and $g'(0)=-w<0$. Moreover, putting $\displaystyle F(\phi)=\int_0^{\phi}g(\xi)d\xi$ and $\phi_0=inf\{\xi>0\,:\, F(\xi)=0\}$, $\phi_0>0$
and $g(\phi_0)=F'(\phi_0)>0$. By applying Theorem 5 and Remark 6.3 in \cite{BerLions}, there exists a unique solution $\phi\in C^3(\er)$ of (\ref{nova21}) such that 
$\phi(0)=\phi_0$, $\phi$ positive and radially decreasing, and such that 
\begin{equation}
\label{equacao3} 
\phi(x),|\phi'(x)|,|\phi''(x)|\leq Ce^{-\delta |x|},
\end{equation}
where $C$ and $\delta$ are positive constants.\\
\\
We can easily deduce from (\ref{equacao1}), (\ref{equacao2}) and (\ref{equacao3}) that $\psi=-\frac 1{\gamma}\phi^{\frac 23}\in C^2(\er)$ with
\begin{equation}
|\psi(x)|,|\psi'(x)|,|\psi''(x)|\leq C'e^{-\frac{2\delta}3 |x|},\quad C'>0.
\end{equation}
\section{Linearized Stability for $p>-\frac 23$}
\noindent
In this section we will consider, for $\displaystyle p>-\frac 23$, special solutions $(\tilde{u},\tilde{v})$ of system  (\ref{sistemainicial}), of the form
\begin{equation}
 \label{especial}
 \left\{\begin{array}{lllll}
         \tilde{u}(x,t)=e^{iwt}e^{i\frac c2(x-ct)}\phi(x-ct)\\
         \tilde{v}(x,t)=\psi(x-ct),
        \end{array}\right.
\end{equation}
satisfying the following conditions:
\begin{itemize}
\item $c\geq 0$ and $c=0$ if $\displaystyle -\frac 23<p<0$; 
\item $\phi,\psi\in C^2(\er)\cap W^{2,\infty}(\er)$;
\begin{equation}
 \label{condicoes}
\end{equation}
\item $\phi,-\psi\geq 0$, and $\phi,-\psi$ radially decreasing;
\item $\phi^p\in C^2(\er)\cap W^{2,\infty}(\er)$ if $p\geq 0$ (cf. Remark \ref{remu}).\\
\end{itemize}
By linearizing the system (\ref{sistemainicial}) around $(\tilde{u},\tilde{v})$ (cf. \cite{Amorim1},\cite{GR}), identifying the first order terms and, for sake of 
simplicity, replacing the solution $(U,V)$ by the new dependent variables $u(x,t)=e^{-iwt}e^{-i\frac{c^2}2t}U(x,t)$ and $v(x,t)=V(x,t)$, we obtain the system
\begin{equation}
 \label{24}
 \left\{\begin{array}{llllllll}
iu_t+u_{xx}&=&(w-\frac{c^2}{2})u+\frac a2\phi^p[(p+2)u+pe^{icx}\overline{u}]+e^{i\frac c2 x}\phi v+\psi u\\
\\
v_t-3\gamma(\psi^2v)_x&=&2Re(e^{i\frac c2 x}\phi u)_x,
        \end{array}\right.
\end{equation}
which we complete with initial data 
\begin{equation}
 \label{id}
 (u_0,v_0)\in H^2(\er)\times H^1(\er).
\end{equation}
Since, for $p<0$, $\phi^p$ is not, in general, a $C^2\cap W^{2,\infty}(\er)$ function, we begin by the study of a regularized system (with the same initial data):
\begin{equation}
 \label{25}
 \left\{\begin{array}{llllllll}
iu_t+u_{xx}&=&(w-\frac{c^2}{2})u+\frac a2(\phi+\epsilon)^p[(p+2)u+pe^{icx}\overline{u}]+e^{i\frac c2 x}\phi v+\psi u\\
\\
v_t-3\gamma(\psi^2v)_x&=&2Re(e^{-i\frac c2 x}\phi u)_x,
        \end{array}\right.
\end{equation}
where $\epsilon>0$ if $p<0$ ($\epsilon=0$ otherwise).\\
\\
We begin by proving the following result concerning this regularized system:
\begin{Propriedade}
 For each $p>-\frac 23$ there exists a unique solution 
 $$(u,v)\in (C([0,+\infty[;H^2)\cap C^1([0,+\infty[;L^2))\times(C([0,+\infty[;H^1)\cap C^1([0,+\infty[;L^2))$$
 of system (\ref{25}) with initial data $(u_0,v_0)\in H^2(\er)\times H^1(\er)$.
\end{Propriedade}
{\bf Proof:}\\
\\
We follow the technique in \cite{DFO},\cite{O} and introduce an auxiliary system with non-local source which can be tackled by Kato's theory (\cite{Kato1}, \cite{Kato2}).
This is necessary in order to write the system (\ref{25}) without derivative loss in the nonlinear term (see \cite{DFO} for details). Hence, we consider the system

\begin{equation}
 \label{26}
 \left\{\begin{array}{llllllll}
 iF_t+F_{xx}&=&(w-\frac{c^2}{2}+\psi)F+a(\phi+\epsilon)^p[\frac p2(F+e^{icx}\overline{F})+F]\\

&&+e^{i\frac c2 x}\phi[3\gamma(\psi^2v)_x+2Re(e^{-i\frac c2x}\phi\tilde{u})_x]+\psi_tu\\
&&+e^{icx}\phi_tv+ap(\phi+\epsilon)^{p-1}\phi_t[\frac p2(u+e^{icx}\overline{u})]\\
v_t-3\gamma(\psi^2v)_x&=&2Re(e^{-i\frac c2 x}\phi \tilde{u})_x,
        \end{array}\right.
\end{equation}
where 
\begin{equation}
 \label{27}
 \left\{\begin{array}{llllllll}
 u(x,t)&=&u_0(x)+\int_0^t F(x,s)ds,\\
\tilde{u}(x,t)&=&(\Delta-1)^{-1}([(w-\frac{c^2}2)+\psi+\frac{a(p+2)}2(\phi+\epsilon)^p]u\\
&&+\frac{ap}2\phi^pe^{icx}\overline{u}+e^{i\frac c2x}\phi v-iF),
       \end{array}\right.
\end{equation}
with initial data
\begin{equation}
 \label{28}
 F(.,0)=F_0\in L^2(\er),\quad v(.,0)=v_0\in H^1(\er).
\end{equation}
Once we have, for a fixed $T>0$, a solution
\begin{equation}
 \label{29}
 F\in C([0,T]; L^2)\cap C^1([0,T];H^{-2}),\quad v\in C([0,T]; H^1)\cap C^1([0,T];L^2)
\end{equation}
for the problem (\ref{26})-(\ref{27})-(\ref{28}), we can argue as in \cite{DFO}, Lemma 2.1, and show that $(u,v)$ is the desired solution to system (\ref{25}). We only
sketch the argument, since it is similar to the one in \cite{Amorim1} and \cite{DFO}.\\
First, we write (\ref{26}) as a system of three equations, by decomposing $F$ into its real and imaginary parts. This allows us to obtain a system with the abstract form
\begin{equation}
 \label{30}
 U_t+AU=g(t,U), \,U(.,0)=U_0,
\end{equation}
with $U=(Re F, Im F, v)$ and $U_0=(Re F_0, Im F_0, v_0)$, the corresponding initial data.\\
Following  \cite{Amorim1}, \cite{DFO}, we decompose the operator
$$A=\left[\begin{array}{cccccccccc}
           0&\Delta&0\\
           -\Delta&0&0\\
           0&0&-3\gamma[(\psi^2)_x+\psi^2\frac{\partial}{\partial x}]
          \end{array}\right]$$
          in the form $SAS^{-1}=A+B$ for some operator $B$. In the present setting, we can choose
$$S=\left[\begin{array}{cccccccccc}
           1-\Delta&0&0\\
           0&1-\Delta&0\\
           0&0&(1-\Delta)^{\frac 12}
          \end{array}\right]$$
Note that $S\,:\,Y=L^2\times L^2\times H^1\to X=H^{-2}\times H^{-2}\times L^2$ is an isomorphism. The relevant properties of $S$ (in particular the ones concerning the
entry $(1-\Delta)^{\frac 12}$) can be found in \cite{Kato1}, Section 8. Observe that the right-hand-side of (\ref{30}) is linear in $U$, hence it is straightforward
to derive the necessary estimates for the source term $g$ and we may finally apply Theorem 2 in \cite{Kato2} (or Theorem 7.1 in \cite{Kato1}) and conclude with the
existence of a unique pair $(F,v)$ satisfying (\ref{26})-(\ref{27})-(\ref{28}), which achieves the sketch of the proof.\hfill$\blacksquare$\\
\\
We are now in position to prove the linearized stability result:\\
\begin{Propriedade}
\label{estavel}
Let $p>-\frac 23$ and consider a special solution $(\tilde{u},\tilde{v})$ to (\ref{sistemainicial}) satisfying (\ref{especial})-(\ref{condicoes}).Then $(\tilde{u},\tilde{v})$
is linearly stable in the sense that for any $T>0$ and any initial data $(u_0,v_0)\in H^1\times L^2$, the system (\ref{24}) admits a unique weak 
solution $(u,v)\in L^{\infty}(0,T;H^1\times L^2)$ such that
\begin{equation}
 \label{31}
\|(u,v)\|^2_{L^{\infty}(0,T;H^1\times L^2)}\leq G_T(\|(u_0,v_0)\|^2_{H^1\times L^2}),
\end{equation}
where $G_{T}\,:\,\er^+\to\er^+$ is a continuous function vanishing at the origin.\\
Moreover, if $(u_0,v_0)\in H^2\times H^1$ and $p\geq 0$, $(u,v)$ is a strong solution satisfying
$$(u,v)\in [C([0,T];H^2)\cap C^1([0,T];L^2)]\times [C([0,T];H^1)\cap C^1([0,T];L^2)]$$
and
\begin{equation}
\label{32}
\|(u,v)\|^2_{L^{\infty}(0,T;H^2\times H^1)}\leq G_T(\|(u_0,v_0)\|^2_{H^2\times H^1}).
\end{equation}
\end{Propriedade}
\noindent
{\bf Proof:}\\
\\
We consider, for fixed $\epsilon$, the solution $(u_{\epsilon},v_{\epsilon})$ of system (\ref{25}) with initial data $(u_{0\epsilon},v_{0\epsilon})\in H^2\times H^1,$
with 
$$(u_{0\epsilon},v_{0\epsilon})\to (u_0,v_0)\quad \textrm{in }H^1\times L^2.$$
In what follows, for simplicity, we will drop the subscript $\epsilon$. By multiplying the first equation in (\ref{25}) by $\overline{u}$ (respectively by $\overline{u_t}$),
taking the imaginary part (respectively the real part) and integrating, we get
$$\frac 12\frac d{dt}\int |u|^2dx=\frac{ap}2 Im\int (\phi+\epsilon)^pe^{i\frac c2x}\overline{u}^2dx+Im\int e^{i\frac c2x}\phi v\overline{u}dx$$
and
$$\frac d{dt}\left\{\frac 12\int |u_x|^2dx+\frac12\left(w-\frac{c^2}2\right)\int |u|^2dx+\frac {a(p+2)}4\int (\phi+\epsilon)^p|u|^2dx+\frac 12 \int \psi|u|^2dx \right\}+$$
$$\frac{ap}2\int (\phi+\epsilon)^pRe\left(e^{icx}\overline{u}\frac{\partial \overline{u}}{\partial t}\right)dx+\int\phi Re\left(e^{i\frac c2x}v\frac{\partial \overline{u}}{\partial t}\right)dx.$$
We have
$$(\phi+\epsilon)^pRe\left(e^{icx}\overline{u}\frac{\partial \overline{u}}{\partial t}\right)=\frac 12(\phi+\epsilon)^p\frac \partial{\partial t}Re\left(e^{icx}\overline{u}^2\right)$$
$$=\frac 12\frac d{dt}\left\{(\phi+\epsilon)^pRe\left(e^{icx}\overline{u}^2\right)\right\}+pc(\phi+\epsilon)^{p-1}\phi'Re\left(e^{icx}\overline{u}^2\right)$$
(recall that $c=0$ if $-\frac 23<p<0$ and $\epsilon=0$ if $p\geq 0$),
$$\phi Re\left(e^{i\frac c2x}v\frac{d\overline{u}}{dt}\right)=\frac{\partial}{\partial t}\left\{\phi Re\left(e^{i\frac c2x}v\overline{u}\right)\right\}+c\phi' Re\left(e^{i\frac c2x}v\overline{u}\right)
- Re\left(e^{i\frac c2x}\phi\overline{u}\frac{\partial v}{\partial t}\right),$$
and, by the second equation in (\ref{24}),
$$ Re\left(e^{i\frac c2x}\phi\overline{u}\frac{\partial v}{\partial t}\right)=Re\left(3\gamma e^{i\frac c2x}\phi\overline{u}(\psi^2v)_x\right)
+2Re\left(e^{i\frac c2x}\phi\overline{u}\right)Re\left(e^{-i\frac c2x}\phi u\right)_x,$$
and so 
$$\int Re\left(e^{i\frac c2x}\phi\overline{u}\frac{\partial v}{\partial t}\right)dx=-3\gamma Re\int \left(e^{i\frac c2x}\phi\overline{u}\right)_x\psi^2vdx.$$
Now, we also derive, from the second equation in (\ref{24}), 
$$\frac 12\frac{d}{dt}\int v^2dx-3\gamma\int(\psi^2v)_xvdx=2\int Re(e^{-i\frac c2x}\phi u)_xvdx.$$
Moreover,
$$\int (\psi^2v)_xvdx=-\int \psi^2vv_xdx=\frac 12\int (\psi^2)_xv^2dx.$$
By applying Cauchy-Schwarz and Gronwall inequalities, it is now easy to obtain the following estimate for $t\in[0,T]$ and where $G_{T}\,:\,\er^+\to \er^+$
is a continuous function vanishing at the origin and indepedent of $\epsilon$:
\begin{equation}
 \label{40}
 \|u_{\epsilon}(t)\|_{H^1}^2+\|v_{\epsilon}(t)\|_{L^2}^2\leq G_T(\|u_0\|_{H^1}^2+\|v_0\|_{H^1}^2),\quad t\in[0,T].
\end{equation}
The first part of the Theorem is now an easy consequence of (\ref{40}) and (\ref{25}), since, by (\ref{40}), there exists a subsequence 
of $\{(u_{\epsilon},v_{\epsilon})\}$ (still denoted $\{(u_{\epsilon},v_{\epsilon})\}$) and $(u,v)\in L^{\infty}(0,T;H^1\times L^2)$ such that 
\begin{itemize}
\item $u_{\epsilon}\rightharpoonup u$ in $L^{\infty}(0,T;H^1)$ weak *;
\item $v_{\epsilon}\rightharpoonup v$ in $L^{\infty}(0,T;L^2)$ weak *;
\item $(u,v)$ satisfies (\ref{31}) and $(u_t,v_t)\in L^{\infty}(0,T;H^{-1}\times H^{-1}).$
\end{itemize}
Hence, $u\in C([0,T];L^2)$, $v\in C([0,T];H^{-1})$, $(u(0),v(0))=(u_0,v_0)$ and $(u,v)$ is a weak solution of (\ref{24}). The uniqueness follows from (\ref{31}).\\
\\
In the case $p\geq 0$, we have $\phi^p\in C^2(\er)\cap W^{2,\infty}$ (cf. Remark \ref{remu}), so we do not need to regularize $\phi$: we can solve directly (\ref{24}) for initial data
$(u_0,v_0)\in H^2\times H^1$. In this case we still obtain estimates of $v_x$, $v_t$, $u_t$ and $u_{xx}$ in $L^2$ to prove (\ref{32}).\\
Differentiating the second equation of the system (\ref{24}), multiplying by $v_x$, and after a few integrations by parts, we otain
\begin{equation}
 \label{41}
 \frac 12\frac{d}{dt}\int (v_x)^2dx+\frac {15\gamma}{2}\int (\psi^2)_x(v_x)^2dx=3\int Re\left(e^{-i\frac c2x}\phi u\right)_{xx}v_xdx.
\end{equation}
From (\ref{41}) and the first equation in (\ref{24}) we deduce, with $G_{T}\,:\,\er^+\to\er^+$ a continuous function vanishing at the origin: 
\begin{equation}
 \label{42}
 \|v_x\|_2^2\leq G_T(\|v_0\|_{H^1}^2)(\|u_{xx}\|_2^2+\|u_x\|_2^2),\quad t\in [0,T].
\end{equation}
Now, the first equation and (\ref{40}) gives
\begin{equation}
 \label{43}
 \|u_{xx}\|_2^2\leq \|u_t\|_2^2+G_T(\|(u_0,v_0)\|_{H^1\times L^2}),\quad t\in[0,T].
\end{equation}
Finally, we differentiate with respect to time the first equation of (\ref{24}), multiply by $\overline{u_t}$ and integrate the imaginary part to obtain 
\begin{equation}
\label{44}
\frac{d}{dt}\|u_t\|_2^2\leq C(\|u\|_2^2+\|v\|_2^2+\|u_t\|_2^2+\|v_t\|_2^2),\quad t\in [0,T].
\end{equation}
From the second equation in (\ref{24}) we also derive
\begin{equation}
\label{45}
\|v_t\|_2^2\leq C(\|v\|_2^2+\|v_x\|_2^2+\|u\|_2^2+\|u_x\|_2^2),\quad t\in [0,T].
\end{equation}
Applying Gronwall's inequality to (\ref{44}) and, by (\ref{40}), (\ref{42}), (\ref{43}) and (\ref{45}), we obtain the estimate (\ref{32}).\hfill$\blacksquare$\\
\\
{\bf Acknowledgements} The authors are grateful to Luis Sanchez for many discussions and were partially supported by FCT (Portuguese Foundation for Science and
Technology) through the grant PEst-OE/MAT/UI0209/2011.


\begin{thebibliography}{xxxxxxxx}
\bibitem{bubbly} L. Akhatov and D. Khismatullin, Long-wave-short-wave interaction in bubbly
liquids, J. Appl. Math. Mech. 63 (1999), 917-926.
\bibitem{Anton} S. Antontsev, J.P. Dias, M. Figueira and F. Oliveira, Non-existence of global solutions For a quasilinear Benney system, J. Math. Fluid Mech. 13 (2011), 213-222.
\bibitem{nAD} P. Amorim and J.P.Dias, A nonlinear model describing a short wave-long
wave interaction in a viscoelastic medium, Quart. Appl. Math. 71(2013),417-432.
\bibitem{Amorim1} P. Amorim, J.P. Dias, M. Figueira and Ph. LeFloch, The Linear Stability of Shock Waves for the Nonlinear Schr\"odinger-Inviscid Burgers System, J. Dyn. Diff. Equations,
 25 (2013), no. 1, 49-69.
\bibitem{Amorim2} P. Amorim and M. Figueira, Convergence of numerical schemes for short wave long wave interaction equations,  J. Hyper. Differential Equations 8 (2011), 777-800.
\bibitem{angulo} J.A. Pava and J.F. Montenegro, Orbital stability of solitary wave solutions of an interaction equation of short and long dispersive waves, 
J. Diff. Equations 174 (2001), 181-199.
\bibitem{BD1} P. B\'egout and J. D\'\i az, On a nonlinear Schr\"odinger equation with a localizing effect, C.R. Acad.
Sci. Paris I 342 (2007), 459-463.
\bibitem{BD2} P. B\'egout and J. D\'\i az, Localizing Estimates of the Support of Solutions of some Nonlinear Schr\"odinger Equations - The Stationary Case,
Ann. Inst. Henri Poincar\'e 29 (2012), 35-58.
\bibitem{B1} D.J. Benney, A general theory for interactions between short and long waves, Stud. Appl. Math. 56 (1977) 81-94.
\bibitem{B2} D.J. Benney, Significant interactions between small and large scale surface waves, Stud. Appl. Math. 55 (1976), 93-106.
\bibitem{BerLions} H. Berestycki and P.L. Lions, Nonlinear scalar field equations I, Arch. Ration. Mech. Anal. 82 (1983) 313-375.
\bibitem{Bekiranov} D. Bekiranov, T. Ogawa and G. Ponce, Interaction equations for short and long dispersive
waves, J. Funct. Anal. 158 (1998), 357-388.
\bibitem{Champeaux} S. Champeaux, D. Laveder, T. Passot, P.-L. Sulem: Remarks on the parallel propagation of small-amplitude dispersive
Alfven waves, Nonlinear Process. Geophys. 6 (1999), 169-178.
\bibitem{microwave} K. Bubke, D.C. Hutchings, U. Peschel, F. Lederer, Modulational instabil-
ity in optical-microwave interaction, Phys. Rev. E 66 (2002) 604-609.
\bibitem{DF} J.P. Dias and M. Figueira, Existence of weak solutions for a quasilinear version of Benney
equations,  J. Hyper. Differential Equations 4 (2007), no. 3, 555-563.
\bibitem{nDFF1} J.P. Dias, M. Figueira and H. Frid, Vanishing viscosity with short wave long wave interactions for systems of conservation laws, Arch. Rat. Mech. Anal. 196
(2010), 981-1010.
\bibitem{nDFF2} J.P. Dias, M. Figueira and H. Frid, Vanishing viscosity with short wave long wave interactions for multi-d scalar conservation laws, J. Diff. Eq. 251 (2011), 492-503.
\bibitem{nDF} J.P. Dias and H. Frid, Short wave-long wave interactions for compressible
Navier-Stokes equations, SIAM J. Math. Analysis 43 (2011), 764-787.
\bibitem{DFO} J.P. Dias, M. Figueira and  F. Oliveira, Existence of local strong solutions for a quasilinear Benney system, C.R. Acad.
Sci. Paris I 344 (2007), 493-496.
\bibitem{nDFO} J.P. Dias, M. Figueira and F. Oliveira, On the Cauchy problem describing an
electron-phonon interaction, Chin. Ann. Math. 32 B (2011), 483-496.
\bibitem{Funakoshi} M. Funakoshi and M. Oikawa, The resonant interaction between a long internal gravity
wave and a surface gravity wave packet, J. Phys. Soc. Japan 52 (1983), 1982-1995.
\bibitem{Grimshaw} R. H. J. Grimshaw, The modulation of an internal gravity-wave packet and the resonance
with the mean motion, Stud. Appl. Math. 56 (1977), 241-266.
\bibitem{GR} E. Godlewski and P. Raviart, An introduction to the linearized stability of solutions of nonlinear hyperbolic systems of conservation laws, Lecture Notes,
Lisbon Summer School, Ellipes, Lisbon, 1999.
 \bibitem{GN1} B.V.Sz. Nagy, Uber Integralgleichungen zwischen einer Funktion und ihrer Ableitung. Acta
Sci. Math. 10 (1941), 64-74.
\bibitem{GN2} R. Killip and M. Visan, Nonlinear Schr\"odinger Equations at critical regularity, Clay Math. Proc. 17 (2009).
\bibitem{Laurencot} Ph. Lauren\c{c}ot, On a nonlinear Schr\"odinger equation arising in the theory of water waves,
Nonlinear Anal. 24 (1995), no. 4, 509-527.
\bibitem{Stollmann} V. Liskevich and P. Stollmann, Schr\"odinger operators with singular complex potentials as generators: existence and stability, Semigroup Forum  60 (2000), 337-343.
\bibitem{Karpman} V.I. Karpman, On the dynamics of Sonic-Langmuir solitons, Phys. Scripta 11 (1975), 263-270. 
\bibitem{Kato1} T. Kato, Linear evolution equations of hyperbolic type, J. Fac. Sci. Univ. Tokyo 17 (1970), 241-258.
\bibitem{Kato2} T.Kato, Quasilinear equations of evolution, with applications to partial differential equations, Lecture Notes in Mathematics 48 (1975), Springer, 25-70.
\bibitem{LeMesurier} B. LeMesurier, Dissipation at singularities of the nonlinear Schr\"odinger equation through
limits of regularisations, Phys. D 138 (2000), 334-343.
\bibitem{Ma} Y.-C. Ma, The complete solution of the long-wave-short-wave resonance equations, Stud.
Appl. Math. 59 (1978), 201-221.
\bibitem{O} F. Oliveira, Stability of the solitons for the one-dimensional Zakharov-Rubenchik equation, Phys. D 175 (2003), 220-240.
\bibitem{Proment} D. Proment and M. Onorato, A note on an alternative derivation of the Benney equations for short wave-long wave interactions, Eur. Journal of Mechs B/Fluids 34 (2012), 1-6.
\bibitem{Sulem} C. Sulem and P.-L. Sulem, The Nonlinear Schr\"odinger Equation, Appl. Math. Sci., vol. 139, Springer, New York (1999).

\bibitem{Redekopp}V.D. Djordjevic, L.G. Redekopp, On the two-dimensional packets of capillary-gravity waves, J. Fluid. Mech. 79 (1977) 703-714.

\bibitem{Tsutsumi1} M. Tsutsumi and S. Hatano, Well-posedness of the Cauchy problem for the long wave–short
wave resonance equations, Nonlinear Anal. 22 (1994), 155-171.
\bibitem{Tsutsumi2} M. Tsutsumi and S. Hatano, Well-posedness of the Cauchy problem for Benney’s first equa-
tions of long wave short wave interactions, Funkcial. Ekvac. 37 (1994), 289-316.
\bibitem{Yajima} N. Yajima and M. Oikawa, Formation and Interaction of Sonic-Langmuir Solitons, Progr. of Theor. Physics 56 (1976), 1719-1739.
\end{thebibliography}
\end{document}